\documentclass[10pt, onesided]{article}

\usepackage{amsthm,amsmath,amssymb,upref,amscd}

\usepackage[isolatin]{inputenc}
\usepackage{lmodern}
\usepackage[T1]{fontenc}
\usepackage{url}

\newtheorem{T}{Theorem}
\newtheorem{Cor}[T]{Corollary}
\newtheorem{TL}[T]{Lemma}
\newtheorem{Prop}[T]{Proposition}
\newtheorem{Remark}{Remark}

\newtheorem{Example}{Example}

\newcommand{\jt}{\tilde{J}}
\newcommand{\vt}{\tilde{V}}
\newcommand{\wt}{\tilde{W}}
\newcommand{\zt}{\tilde{Z}}
\newcommand{\at}{\tilde{A}}

\newcommand{\tr}{\mathop{\mathrm{tr}}}
\newcommand{\cof}{\mathop{\mathrm{cof}}}
\newcommand{\R}{\mathbb{R}}
\newcommand{\ed}{\mathrm{d}}
\newcommand{\vektor}[2]{\left[\begin{array}{c} #1 \\ #2
    \end{array}\right]} 
\newcommand{\vektort}[3]{\left[\begin{array}{c} #1 \\ #2 \\ #3
    \end{array}\right]}   
\newcommand{\matris}[4]{\left[\begin{array}{cc} #1 & #2 \\ #3 &
      #4\end{array}\right]} 
\newcommand{\matrist}[9]{\left[\begin{array}{ccc} #1 & #2 & #3 \\ #4 &
      #5 & #6 \\ #7 & #8 & #9\end{array}\right]}
\newcommand{\pd}[2]{\frac{\partial #1}{\partial #2}}

\begin{document} \bibliographystyle{abbrv}

\title{Multiplication of solutions for linear overdetermined systems
  of partial differential equations}   
\author{Jens Jonasson\\ Department of Mathematics\\ Link\"oping
  University\\ SE-581 83 Link\"oping, Sweden\\ jejon@mai.liu.se}   
\maketitle

\begin{abstract}

A large family of linear, usually overdetermined, systems of partial
differential equations that admit a multiplication of solutions, i.e,
a bi-linear and commutative mapping on the solution space, is studied.
This family of PDE's contains the Cauchy--Riemann equations and the
cofactor pair systems, included as special cases. The multiplication
provides a method for generating, in a pure algebraic way, large
classes of non-trivial solutions that can be constructed by forming
convergent power series of trivial solutions.

\end{abstract}

\section{Introduction}

In this paper we study a wide class of linear first order systems of
partial differential equations, that allow a bi--linear multiplication
in the space of solutions. The simplest example is the Cauchy--Riemann
equations. We know that two holomorphic functions, $f=V+\mathrm{i}\vt$
and $g=W+\mathrm{i}\wt,$ can be multiplied in order to produce a new
holomorphic function $fg = VW - \vt\wt + \mathrm{i}(V\wt + \vt W)$. In
terms of the Cauchy--Riemann equations
\begin{align*}
  \left\{ 
    \begin{aligned}
      \pd{V}{x} & = \;\pd{\vt}{y}\\
      \pd{V}{y} & = \;-\pd{\vt}{x}\\
    \end{aligned}
  \right.
\end{align*}
this multiplication can be expressed in the following way: two
solutions $(V,\vt)$ and $(W,\wt)$ prescribe, in a bi--linear way, a
new solution $(VW-\vt\wt, V\wt+\vt W).$ From the basic theory of
holomorphic functions, we know that any solution of the
Cauchy--Riemann equations can be expressed locally as a convergent
power series of a simple solution with respect to the described
multiplication.\newline 

The Cauchy--Riemann equations provide the simplest example of a system
of PDE's that has a multiplication on its solution set, but there are
more sophisticated examples. One such example is the multiplication of
cofactor pair systems, discovered by Lundmark in
\cite{lundmark-2001}.\newline

A cofactor pair system (or bi-cofactor system) is a dynamical system
$\ddot{q}^h+\Gamma^h_{ij}\dot{q}^i\dot{q}^j=F^h,$ $h=1,\ldots ,n,$ on
a (pseudo-) Riemannian manifold, such that the force $F=F(q)$ has two
different cofactor formulations $F(q) = (\cof J)^{-1}\nabla V =
\linebreak (\cof \jt)^{-1}\nabla \vt,$ where $J$ and $\jt$ are
independent special conformal Killing tensors of type $(1,1)$, $V$ and
$\vt$ are smooth real-valued functions, $\cof J=(\det{J})J^{-1},$ and
$\nabla$ is the gradient ($(\nabla V)^i = g^{ij}\partial_j V$).
Cofactor pair systems have several desirable properties, in general
they are completely integrable, they admit a bi-Hamiltonian
formulation, and they are equivalent (or correspondent) to separable
Lagrangian systems \cite{benenti-2005, crampin-2001, lundmark-2003,
lundmark-2002, rauch-1999, rauch-2003}. \newline

A cofactor pair system is characterized by a pair of functions $V$ and
$\vt,$ and a pair of special conformal Killing tensors $J$ and $\jt,$
that satisfy the relation
\begin{align}\label{qcr}
  (\cof J)^{-1}\nabla V = (\cof \jt)^{-1}\nabla \vt.
\end{align}  
For fixed special conformal Killing tensors $J$ and $\jt,$ the
equation (\ref{qcr}) constitutes a system of first order linear PDE's
for two functions $V$ and $\vt.$ In \cite{lundmark-2001}, Lundmark
found that the equation (\ref{qcr}) allows a multiplication of
solutions. When $n=2$ the multiplication formula is given by
\begin{align*} 
  (V,\vt)*(W,\wt) = \left(VW-\det{(\jt^{-1}J)}\vt\wt, \;V\wt+\vt W-\tr
    (\jt^{-1}J)\vt\wt \right),
\end{align*}
where $(V,\vt)$ and $(W,\wt)$ are solutions of (\ref{qcr}). We see
that when $\det{(\jt^{-1}J)}$ and $\tr(\jt^{-1}J)$ are not both
constant, we can choose trivial (constant) solutions $(V,\vt)$ and
$(W,\wt)$ of (\ref{qcr}) and obtain non--trivial solutions through the
multiplication. When $n>2$ a multiplication also exists, but one has
to consider the related parameter--dependent system
\begin{align*}
  \left(\cof (J+\mu\jt)\right)^{-1}\nabla V_\mu = (\cof
  \jt)^{-1}\nabla \vt, 
\end{align*}  
which can also be written as
\begin{align*}
  (\jt^{-1}J+\mu I)\nabla V_\mu = \det{(\jt^{-1}J+\mu I)}\nabla \vt,
\end{align*}  
where $V_\mu$ is polynomial in the real parameter $\mu$ (note that
throughout this paper, we use the notation $V_\mu$ rather than
$V(\mu)$ to indicate dependence on the parameter $\mu$). \newline

The most interesting property of this multiplication is that it
provides a tool for producing new cofactor pair systems from known
ones. Especially, infinite families of separable potentials can be
constructed. For example, the Jacobi, Neumann, and parabolic families
of separable potentials are all constructed in \cite{rauch-1986}
through a recursive process that is a special case of the
multiplication of cofactor pair systems.\newline

\begin{Remark}
  In \cite{jodeit-1990}, equations of the form $($\ref{qcr}$)$,
  considered on a real or complex vector space where $J$ and $\jt$ are
  constant matrices, are studied, and the general analytic solution is
  described . 
\end{Remark}

In order to gain better understanding of this multiplication, systems
of the form
\begin{align}\label{xmueq}
  (X+\mu I)\nabla V_\mu = \det{(X+\mu I)}\nabla \vt,
\end{align}  
defined on a general (pseudo-) Riemannian manifold, were studied in
\cite{jonasson-2006}, without referring to any underlying dynamical
system. By analyzing the corresponding equations at each degree of
$\mu$ in the equation (\ref{xmueq}), it becomes obvious that the
equation is satisfied if and only if the degree of $V_\mu$ is $n$ and
the left hand side $(X+\mu I)\nabla V_\mu$ can be written as a product
of the scalar $\det{(X+\mu I)}$ and some $1-$form which is constant in
$\mu.$ We can therefore rewrite the equation (\ref{xmueq}) as
\begin{align}\label{qcrmu}
  (X+\mu I)\nabla V_\mu \equiv 0\quad \left(mod\;\det{(X+\mu
      I)}\right). 
\end{align}  
It turned out that the system (\ref{xmueq}) allows for a
multiplication of solutions, similar to the one existing for cofactor
pair systems, if and only if the tensor $X$ satisfies the equation
\begin{align}\label{xmueq2}
  (X+\mu I)\nabla \det{(X+\mu I)} = \det{(X+\mu I)}\nabla \tr (X+\mu
  I). 
\end{align} 
Several classes of solutions of (\ref{xmueq2}) where discovered, and
it became apparent that systems of the form (\ref{qcr}) and the
Cauchy--Riemann equations only constitute special cases of a much
larger family of systems of PDE's that admit a multiplicative
structure on the solution space. \newline

It was also remarked in \cite{jonasson-2006}, that by considering more
general systems than (\ref{xmueq}), one finds other new classes of
systems that allow multiplication. In this paper, we will examine that
subject. The linear systems of PDE's that we consider are in general
impossible to solve, but the multiplication provides a non-trivial
superposition principle (on top of the ordinary linear superposition)
that, given to solutions, prescribes a new solution in a bi-linear and
pure algebraic way. With this superposition principle, large classes
of new solutions can be generated from known solutions. In particular,
we can construct non--trivial solutions by forming convergent power
series of a simple solution. The question then arises for which
systems of linear PDE's these power series constitute all solutions,
like in the case of the Cauchy--Riemann equations where all
holomorphic functions admit a power series representation. Besides
providing us with more systems of PDE's that admit a multiplicative
structure on the solution set, the generalization helps us to better
understand the multiplication for the systems already known (in
particular the puzzling multiplication of cofactor pair
systems).\newline

This paper is organized as follows. In section \ref{sec_mult} we
formulate an abstract framework for characterizing the class of
systems of PDE's that admit multiplication. We define the $*-$operator
and give a characterization of those systems that admit
$*-$multiplication on the set of solutions. The multiplication
provides a method for generating, in a pure algebraic and non-trivial
way, new solutions from known solutions. A second formulation of the
systems, using related matrices, is introduced. In this matrix
notation the Euclidean algorithm for polynomial division, which is
closely related to the $*-$multiplication, can be encoded in an
explicit polynomial of matrices. Some algebraic properties of the
multiplication are also mentioned in this section. In section
\ref{sec_exp} we investigate the explicit forms of systems that admit
multiplication. The study splits into different cases depending on
relations among certain discrete parameters that appear in the studied
class of systems. Especially, some typical (generic) systems with
multiplication are derived and examined. The most interesting property
of the multiplication of solutions is that we can construct large
classes of solutions by forming power series, with respect to the
$*-$multiplication, of trivial solutions. Section \ref{sec_ps} is
devoted to study such power series solutions. The problem of
constructing systems with $*-$multiplication is in general quite
complicated. In section \ref{sec_find}, several methods for
constructing systems with $*-$multiplication are described. The last
section \ref{sec_conc} contains concluding remarks and natural
questions raised by the study presented in this paper.

\section{Multiplication of solutions for linear systems of
  PDE's }\label{sec_mult}

Let $Q$ be a $n-$dimensional differentiable real manifold. Consider
equations of the form  
\begin{align}\label{sys}
  A_\mu\ed V_\mu \equiv 0 \quad (mod\;Z_\mu),
\end{align}
where $A_\mu$ is a $(1,1)-$tensor depending polynomially on the real
parameter $\mu,$ and $Z_\mu,$ $V_\mu$ are real-valued $\mathcal{C}^1$
functions on $Q$ that also depend polynomially on $\mu$. The
expression $A_\mu\ed V_\mu$ is a $1-$form which components are
polynomial in $\mu,$ and the unknown function $V_\mu$ is a solution of
the equation (\ref{sys}) if these components are all divisible (when
considered as polynomials in $\mu$) by the fixed function $Z_\mu$. Let
$Z_\mu = Z_0 + \mu Z_1 + \cdots + \mu^{m-1}Z_{m-1} + \mu^m,$ be a
polynomial of degree $m$, then there is no restriction to assume that
$A_\mu = A_0 + \mu A_1 + \cdots + \mu^kA_k$ has degree at most $m-1$
(otherwise we can reduce it modulo $Z_\mu$). In order to simplify the
description of systems admitting a multiplication, we also assume that
$V_\mu = V_0 + \mu V_1 + \cdots + \mu^{m-1}V_{m-1}$ has degree $m-1.$

\begin{Remark}
  The system $($\ref{sys}$)$ should be compared with system
  $($\ref{qcrmu}$)$.  We see that $($\ref{sys}$)$ generalizes
  $($\ref{qcrmu}$)$ in several ways.  First of all, no metric is
  specified on the manifold $Q$ corresponding to the system
  $($\ref{sys}$)$, and we consider a system of equations for $1-$forms
  using the exterior differential operator $\ed,$ rather than a system
  of equations for vector fields expressed with the gradient operator
  $\nabla.$ Moreover, in $($\ref{qcrmu}$)$ the polynomial $X + \mu I$
  is of degree one and contains only one arbitrary tensor $X$, while in
  $($\ref{sys}$)$ we consider a more general polynomial $A_\mu$, and
  instead of $\det{(X + \mu I)},$ we consider an arbitrary polynomial
  $Z_\mu$ that does not have to be related to $A_\mu.$ 
\end{Remark}

Since the highest order coefficient of $Z_\mu$ is a unit, for each
polynomial $P_\mu$ with coefficients in the commutative ring of
real--valued functions on $Q,$ there exists unique polynomials $Q_\mu$
and $R_\mu$ such that $P_\mu = Q_\mu Z_\mu + R_\mu,$ where
$\deg{R_\mu}<m.$ Thus, for each function $V_\mu,$ there exists unique
$1-$forms $B_0, \ldots ,B_{m-1}$ such that
\begin{align}\label{B_sys}
  A_\mu\ed V_\mu \equiv B_0 + \mu B_1 + \cdots + \mu^{m-1}B_{m-1}
  \quad (mod\;Z_\mu),   
\end{align}
and the equation (\ref{sys}) can be written as $B_0 = B_1 =
\cdots = B_{m-1} = 0$. Thus, in local coordinates $q^1,\ldots,q^n,$
the system (\ref{sys}) constitutes a, usually overdetermined, system
of $nm$ first order linear partial differential equations for $m$
dependent variables. We will see that there exist non-trivial
systems (\ref{sys}) admitting a multiplicative structure on its
solution set.\newline

Define a bilinear operation $*,$ on the set of all real-valued
functions on $Q$ that are polynomial in $\mu,$ by letting
$V_\mu*W_\mu$ be the residue of the ordinary product $V_\mu W_\mu$
modulo $Z_\mu.$ In other words, $V_\mu*W_\mu$ is the unique polynomial
of degree less than $m$ that can be written as $V_\mu W_\mu - Q_\mu
Z_\mu,$ for some polynomial $Q_\mu.$ For certain choices of $A_\mu$
and $Z_\mu$, the $*-$multiplication maps solutions of (\ref{sys}) to
new solutions:

\begin{T}[$*-$Multiplication]\label{mult_thm}
   Let $S$ denote the solution set of $($\ref{sys}$)$. Then $*$ is a
   bilinear operation on $S$ if and only if 
   \begin{align}\label{chareq}
     A_\mu\ed Z_\mu\equiv 0 \quad (mod\;Z_\mu),
   \end{align}
   i.e., if and only if $Z_\mu - \mu^n\in S$.
\end{T}

\begin{proof}
  Given two solutions $V_\mu, W_\mu\in S$, let $Q_\mu$ be the
  polynomial such that the product $V_\mu W_\mu$ can be written as
  $V_\mu W_\mu = Q_\mu Z_\mu + V_\mu*W_\mu.$ Then, we have
  \begin{align*}
    A_\mu\ed (V_\mu*W_\mu) & = W_\mu A_\mu\ed V_\mu + V_\mu A_\mu\ed
    W_\mu - Z_\mu A_\mu\ed Q_\mu - Q_\mu A_\mu\ed Z_\mu \\ 
    & \equiv -Q_\mu A_\mu\ed Z_\mu \quad (mod\;Z_\mu).
  \end{align*}
  Thus, we see that $A_\mu\ed Z_\mu\equiv 0$ is a sufficient condition
  for the existence of the bi--linear operation $*$ on $S$. To see that
  it is also a necessary condition, consider the trivial solutions
  $V_\mu=\mu$ and $W_\mu=\mu^{m-1}.$ For this choice of solutions, the
  polynomial $Q_\mu$ becomes a non--zero constant, which forces the
  relation (\ref{chareq}) to be satisfied in order for $V_\mu*W_\mu$
  to be a solution.
\end{proof}

The following algebraic properties of $*$ are immediate consequences
of the corresponding properties of multiplication in general quotient
rings of polynomials: 

\begin{Cor}
  The solution set $S$ together with the scalar multiplication,
  addition (defined in the obvious way) and multiplication $*$ is an
  algebra over $\R$, where the $*$--multiplication is associative and
  commutative.
\end{Cor}

To calculate the $*-$product $V_\mu*W_\mu,$ we form the ordinary
product  
\begin{align*}
  V_\mu W_\mu & = V_0W_0 + \left(V_0W_1 + V_1W_0\right)\mu + \cdots +
  \\ 
  & \quad \left(V_{m-2}W_{m-1}+V_{m-1}W_{m-2}\right)\mu^{2m-3} +
  V_{m-1}W_{m-1}\mu^{2m-2},   
\end{align*}
and replace $\mu^m,\ldots , \mu^{2m-2}$ with their residues modulo
$Z_\mu:$
\begin{align*}
  \mu^m & \equiv -Z_0 - Z_1\mu - \cdots - Z_{m-1}\mu^{m-1}\\ 
  \mu^{m+1}& \equiv Z_0Z_{m-1} + \left(Z_1Z_{m-1}-Z_0\right)\mu +
  \cdots + \left(Z_{m-1}Z_{m-1}-Z_{m-2}\right)\mu^{m-1} \\ 
  &\;\; \vdots
\end{align*}
In general it is hard to find non--trivial solutions of the system
(\ref{sys}), but having the $*-$operator we can generate (in a pure
algebraic way) an infinite family of non-trivial solutions by starting
with trivial solutions. For example, we can construct non--trivial
solutions by forming polynomials, or convergent power series, of the
trivial solution $\mu\in S:$
\begin{align}\label{ps}
  V_\mu = \sum_r a_r\mu^r_*, \quad \mathrm{where} \quad \mu^r_* :=
  \underbrace{\mu*\cdots *\mu}_{r\textrm{ factors}}, 
\end{align}
with $a_r$ being real constants. We note that, since all trivial
solutions can be expressed as polynomials in $\mu,$ every solution
that is a sum of products of trivial solutions has the form (\ref{ps})
again. As long as we consider domains in the manifold $Q$ where
$Z_0\ne 0,$ we can also allow negative powers $\mu^{-a}_*$ in
(\ref{ps}), for every natural number $a$, by defining $\mu^{-a}_* :=
(\mu^{-1})^a_*$ and
\begin{align*}
  \mu^{-1} := -\frac{1}{Z_0}\left(Z_1 + \cdots +
    Z_{m-1}\mu^{m-2}+\mu^{m-1}\right),  
\end{align*}
so that $\mu*\mu^{-1} = 1$.\newline

The following two examples illustrate how the Cauchy--Riemann
equations and the cofactor pair systems can be considered as special
cases of systems of the form (\ref{sys}) that admit a
$*-$multiplication of the kind described in theorem \ref{mult_thm}:

\begin{Example}[Cauchy--Riemann equations]
  Let $Q$ be the $2-$dimensional Euclidean space with Cartesian
  coordinates $(x,y)$, $Z_\mu = 1+\mu^2$ and  
  \begin{align*}
    A_\mu = \matris{0}{-1}{1}{0} + \mu I,
  \end{align*}
  where $I$ is the identity matrix. Then the system $($\ref{sys}$)$
  reduces to the Cauchy--Riemann equations, and the ordinary
  multiplication of holomorphic functions follows from the
  $*-$multiplication:
  \begin{align*}
    (V,\vt)*(W,\wt) = \left(VW - \vt\wt, \;V\wt+\vt W \right).    
  \end{align*}
 Every solution $(V,\vt)$ that satisfies the Cauchy--Riemann equations
 in the origin, can be expressed in a neighborhood of the origin as a
 power series of the simple solution $(x,y)$ :
  \begin{align*}
    V_{\mu} = V + \mu\vt = \sum_{r=0}^\infty a_r (x + y\mu )_*^r,   
  \end{align*}
  where $(x + y\mu )_*^r$ again denotes the $r$'th power with respect
  to the $*-$multiplication.
\end{Example}

\begin{Example}[Multiplication of cofactor pair systems]  
  Let $m=n$ and suppose that the tensor $A_\mu$ is linear in $\mu$ and
  has the identity mapping as the highest order coefficient, i.e.
  $A_\mu = X + \mu I.$ The system $($\ref{sys}$)$ can for this special
  case be written as 
  \begin{align}\label{oldsys}
    X\ed V_i + \ed V_{i-1} = Z_i\ed V_{n-1},\quad i = 0,\ldots ,n-1,
  \end{align}
  where $V_{-1}:=0$. If we also let $Z_\mu = \det{(X + \mu I)},$ the
  system $($\ref{sys}$)$ reduces to the system $($\ref{xmueq}$)$ when
  we specify a metric on $Q$ and consider the equivalent ``vector
  version'' $A_\mu\nabla V_\mu\equiv 0$ of $($\ref{sys}$)$.
  Restricting the attention to the case $Z_\mu = \det{(X + \mu I)}$ is
  quite natural since, if we also assume that the coefficients
  $Z_0,\ldots ,Z_{n-1}$ of $Z_\mu$ are functionally independent, it is
  a necessary condition for the equation $($\ref{chareq}$)$ to be
  satisfied. To see this, choose coordinates $q^1,\ldots ,q^n$ as $q^i
  = Z_{i-1}.$ Then the equation $($\ref{chareq}$),$ or equivalently
  $X\ed Z_i + \ed Z_{i-1} = Z_i\ed Z_{n-1},$ reduces to $X^i_j =
  q^i\delta^n_j - \delta^{i-1}_j,$ where $\delta$ is the Kronecker
  delta symbol. In other words, $-X$ must in these coordinates be the
  companion matrix (see $($\ref{compmatrix}$)$) of the polynomial
  $Z_\mu,$ and therefore it follows that $Z_\mu = \det{(X + \mu
    I)}.$\newline 
  
  The $*-$multiplication reduces to Lundmarks multiplication of
  cofactor pair systems if we let $X=\jt^{-1}J,$ where $J$ and $\jt$ are
  special conformal Killing tensors. In \cite{jonasson-2006}, several
  other families of tensors $X,$ that satisfy the equation
  $($\ref{xmueq2}$)$ have been found.
\end{Example}

\subsection{Matrix notation}\label{matrixnot}

For the purpose of further study of systems of PDE's of the form
(\ref{sys}), we shall introduce a new kind of matrix formulation for
these systems and for the corresponding $*-$multiplication. The matrix
formulation makes it possible to give an explicit formula for
calculating powers of solutions, with respect to the
$*-$multiplication. \newline

The idea is to consider the column matrix $V = [V_0, V_1, \ldots,
V_{m-1}]^T$ instead of the polynomial $V_\mu = V_0 + \mu V_1 + \cdots
+ \mu^{m-1}V_{m-1},$ and to observe that
\begin{align*}
  V = V_Ce_1 = \left( V_0C^0 + V_1C^1 + \cdots + V_{m-1}C^{m-1}
  \right) e_1 
\end{align*}
with $V_C = V_{\mu=C}$ where we have formally substituted the
parameter $\mu$ with the companion matrix
\begin{align}\label{compmatrix}
  C := C[Z_\mu] = \left[
  \begin{array}{ccccc}
    0 & 0 & \cdots & 0 & -Z_0 \\
    1 & 0 & \cdots & 0 & -Z_1 \\
    0 & 1 &  \ddots & & \vdots \\
    \vdots & \vdots & \ddots & 0 & -Z_{n-2} \\
    0 & 0 & \cdots & 1 & -Z_{n-1} 
  \end{array}
  \right]
\end{align} 
of $Z_\mu$, and where $e_1 = [1,0,\ldots, 0]^T.$ The most important
advantage of the matrix notation is that we can express the Euclidean
algorithm of polynomial division in a more explicit way. For any
polynomial $P_\mu = P_0 + \cdots + P_t\mu^t,$ the residue modulo
$Z_\mu$ can be written as $R_\mu = R_0 + \cdots + R_{m-1}\mu^{m-1},$
where $[R_0, R_1, \ldots, R_{m-1}]^T = P_Ce_1.$ Thus, in the matrix
notation, the $*-$multiplication of two solutions $V$ and $W$ of
(\ref{sys}) can be written as $V*W = V_CW_Ce_1.$ In particular, the
matrix version of $\mu^a_*$ is $C^ae_1.$ The system (\ref{sys}) can
also be expressed in terms of matrices as
\begin{align}
  0 & =
  \sum_{i=0}^k\sum_{j=0}^{m-1}C^{i+j}e_1\left[\partial_1V_j,\ldots, 
    \partial_nV_j\right]A_i \nonumber\\ 
  & = \sum_{i=0}^k C^i \left(\sum_{j=0}^{m-1}
    C^je_1\left[\partial_1V_j, \cdots, \partial_nV_j\right]
  \right)A_i \nonumber\\  
  & = \sum_{i=0}^k C^iV'A_i, \label{matrix_sys}
\end{align}
where we consider $A_i$ as the matrix with elements
$(A_i)^a_b$ and $V'$ is the functional matrix 
\begin{align*}
  V' =
  \matrist{\partial_1V_0}{\cdots}{\partial_n
    V_0}{\vdots}{}{\vdots}{\partial_1
    V_{m-1}}{\cdots}{\partial_nV_{m-1}}.   
\end{align*} 
The equation (\ref{matrix_sys}) is indeed independent of coordinates.
The expression on the right-hand side of equation (\ref{matrix_sys})
is a $m\times 1-$matrix consisting of the $1-$forms $\sum_{i,j}
(C^i)_{aj}A_i\ed V_{j-1},$ $a = 1,\ldots ,m,$ where $(C^i)_{aj}$ is
the element in the row $a$ and the column $j$ of the matrix $C^i.$ The
$*-$multiplication theorem can then be expressed as:

\begin{Prop}[Explicit criterion for existence of $*-$multiplication]
  Let $V$ and $W$ be solutions of $0 = \sum_{i=0}^k C^iV'A_i,$
  then $V*W := V_CW_Ce_1$ is also a solution if and only if 
  \begin{align}\label{chareqmatrix}
    0 = \sum_{i=0}^kC^iZ'A_i.
  \end{align}
\end{Prop}

\section{Explicit form of linear PDE's admitting\\
  $*-$multiplication}\label{sec_exp}

The system (\ref{sys}) contains three parameters: $n$ -- the dimension
of the manifold $Q;$ $m$ -- the polynomial degree of the function
$Z_\mu;$ $k$ -- the polynomial degree of the tensor $A_\mu.$ In this
section we will discuss how the form of the system (\ref{sys}), or
(\ref{matrix_sys}) in matrix notation, and of the related
$*-$multiplication depends on these numbers. We will also, for
different choices of $n,m,k,$ specify the structure of typical
(generic) systems that allow $*-$multiplication. This is done by
choosing the functions $Z_i$ as coordinates.\newline

When $k=0,$ the system (\ref{sys}) can be written as $A_0\ed V_i=0,$
$i=0,\ldots,m-1$ and the multiplication becomes trivial. When $n=1,$
the system (\ref{sys}) reduces to a quite simple system of ordinary
differential equations. Therefore, we consider only cases where $k>0$
and $n>1.$\newline

The $*-$product $V*W$ of two solutions is a collection of $m$
functions, each being a sum of functions of the form $P_{ij}V_iW_j,$
where $P_{ij}$ is a polynomial expression of the variables $Z_k.$ The
structure of the $*-$multiplication formula depends only on the
parameter $m,$ not on $n$ or $k.$ Since the degree of the polynomials
$P_{ij}$ will not exceed $m-1,$ the multiplication will be more
complex for higher values of the parameter $m.$ For the simplest case
$m=2,$ we have
\begin{align*}
  V*W & = V_CW_Ce_1 \\
  & = \left( V_0W_0C^0 + \left(V_0W_1 + V_1W_0\right)C + V_1W_1C^2
  \right)e_1 \\ 
  & = \vektor{ V_0W_0 -  Z_0V_1W_1}{ V_0W_1 + V_1W_0 -  Z_1V_1W_1}, 
\end{align*}
and for $m=3$ the $*-$product $V*W$ is given by   
\begin{align}\label{mult3}
  \vektort{V_0W_0 - Z_0V_1W_2 - Z_0V_2W_1 + Z_0Z_2V_2W_2}{
    V_0W_1 + V_1W_0 - Z_1V_1W_2 - Z_1V_2W_1 + \left(-Z_0 +
      Z_1Z_2\right)V_2W_2}{ V_0W_2 + V_1W_1 + V_2W_0 - Z_2V_1W_2 -
    Z_2V_2W_1 + \left(-Z_1 + Z_2^2\right)V_2W_2} 
\end{align}

\subsection{Generic cases for different choices of $(n,m,k)$}

Our approach to find explicit forms of equations admitting
$*-$multiplication, is to choose some generic coordinates in which
the system (\ref{sys}), equipped with $*-$multiplication, takes a
simple form. Since the functions $Z_i$ play a fundamental role for
the multiplication, we will assume that as many of these functions as
possible are functionally independent and take them as local
coordinates on $Q$. Since there are at most $n$ different functionally
independent functions, the relation between $m$ (the number of
functions $Z_i$) and $n$ will be crucial for specifying each
generic case.\newline

When $m=n,$ we can choose as generic coordinates $q^1 = Z_0,\ldots
,q^n = Z_{m-1}$ if the functions $Z_i$ are functionally independent.
We consider the case when $m=n$ as our main case since the system
(\ref{sys}) takes a simpler form than in the other cases.\newline

When $m<n,$ the functions $Z_i$ are too few to form a complete set of
coordinates. Instead we choose generic coordinates $q^1, \ldots , q^n$
such that $q^1 = Z_0,\ldots ,q^m = Z_{m-1},$ without specifying the
last $n-m$ coordinates $q^{m+1}, \ldots , q^n.$\newline

When $m>n,$ the functions $Z_i$ must be functionally dependent. For
the generic case we assume that $Z_0, \ldots , Z_{n-1}$ are
functionally independent, and choose them as generic
coordinates.\newline

We shall present below an explicit form of the system (\ref{sys}) in
generic coordinates for the cases $m=n,$ $m<n,$ and $m>n.$ For each
case we will also consider simpler sub-cases according to the
following schematic diagram (\ref{subcases}) for the triples $(n,m,k)$
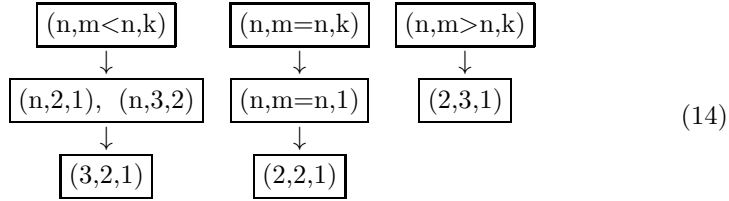
\begin{align}\label{subcases}
  \begin{array}{ccc}
    \fbox{(n,m<n,k)} & \fbox{(n,m=n,k)} & \fbox{(n,m>n,k)}\\
    \downarrow & \downarrow & \downarrow \\
    \fbox{(n,2,1),\; (n,3,2)} & \fbox{(n,m=n,1)} & \fbox{(2,3,1)}\\
    \downarrow & \downarrow &  \\
    \fbox{(3,2,1)} & \fbox{(2,2,1)} &\\\\
  \end{array}
\end{align}

For the sake of simplicity, we will consider the case when $A_{m-1}$
is non--singular. It is then no restriction to assume that $A_{m-1}$
is the identity mapping.

\subsection{$m=n$}

A generic case with the simplest structure is obtained when $m = n.$
According to the discussion above, $(n,m,k) = (2,2,1)$ is the lowest
value of the parameters for which the multiplication is non-trivial.
It is also the best case to study in order to get a good understanding
of the mechanism of the multiplication. We will investigate this case
in detail, and after that some of the ideas will be generalized to the
cases $(n = m,m,k)$ and $(n = m,m,k =1)$. 

\subsubsection{$(n,m,k) = (2,2,1)$}

For this choice of parameters $n,$ $m$ and $k$, we have $Z_\mu = Z_0 +
\mu Z_1 + \mu^2,$ $A_\mu = A_0 + \mu A_1,$ $V_\mu = V_0 + \mu V_1$ and
(\ref{sys}) can be written as
\begin{align}\label{sys221}
  \left\{
    \begin{aligned}
      A_0\ed V_0 & = \; Z_0A_1\ed V_1 \\
      A_1\ed V_0 & = \; \left( Z_1A_1 - A_0
      \right) \ed V_1,
    \end{aligned}
  \right.
\end{align}
or, as we have seen, in matrix notation as $0 = V'A_0 + CV'A_1.$ Thus,
in local coordinates, (\ref{sys221}) constitutes a system of four
partial differential equations for two unknown functions $V_0,$ $V_1$
of two independent variables $x,y.$ Since the number of equations
exceeds the number of dependent variables, this system will in general
be overdetermined. If we assume that $A_1$ is non-singular, we can
instead of (\ref{sys221}) consider the equivalent system
\begin{align}\label{sys221b}
  \left\{
    \begin{aligned}
      0 & = \; \left( A^2 - Z_1 A + Z_0 I \right)\ed V_1 \\
      \ed V_0 & = \; \left( Z_1 I - A \right) \ed V_1,
    \end{aligned}
  \right.
\end{align}
where $A = A_0A_1^{-1},$ or in matrix notation $0 = V'A + CV'.$ We
assume now that the functions $x = Z_0,$ $y = Z_1$ are functionally
independent and consider the system (\ref{sys221b}) in the generic
coordinates $x,y.$ In these coordinates the relation
(\ref{chareqmatrix}), that guarantees a $*-$multiplication of the
corresponding system (\ref{sys221b}), reduces to $A = -C.$ Thus, in
generic coordinates, if we require existence of $*-$multiplication,
the first equation of system (\ref{sys221b}) is a consequence of the
Cayley--Hamilton theorem so the system reduces to the system $\ed V_0
= \left( y I + C \right) \ed V_1,$ that has components:
\begin{align*}
  \left\{
    \begin{aligned}
      \pd{V_0}{x} & = y\pd{V_1}{x} + \pd{V_1}{y}\\ 
      \pd{V_0}{y} & = -x\pd{V_1}{x}.
    \end{aligned}
  \right.
\end{align*} 
Thus, the generic $(2,2,1)-$case constitutes in fact a determined
system of two partial differential equations for two unknown functions
of two independent variables. As we have seen, the $*-$product of two
solutions can in this case be written as
\begin{align*}
  V*W = V_CW_Ce_1 = \vektor{ V_0W_0 -  xV_1W_1}{ V_0W_1 + V_1W_0 -
    yV_1W_1}. 
\end{align*}
The simplest non-trivial solutions obtained by taking powers of the
trivial solution $(0,1)$ are
\begin{align}\label{tp}
  \begin{aligned}
    (0,1)^{-2}_* & = \left( -\frac{y^2-x}{x^2},
      -\frac{y}{x^2} \right)\\  
    (0,1)^{-1}_* & = \left( -\frac{y}{x}, -\frac{1}{x} \right)\\ 
    (0,1)^0_* & = (1,0)\\
    (0,1)^1_* & = (0,1)\\
    (0,1)^2_* & = \left( -x,-y \right)\\
    (0,1)^3_* & = \left( xy, -x + y^2 \right)
  \end{aligned}
\end{align}
If the roots of the polynomial $Z_\mu$ (or the eigenvalues of the
companion matrix $C$) are functionally independent, we can instead
define local coordinates through $Z_0 = xy,$ $Z_1 = x + y.$ The
condition (\ref{chareq}) can then be expressed as $A_0 = DA_1,$ where
$D = \mathrm{diag}(x,y).$ If we moreover assume that $A_1$ is
non-singular, the system (\ref{sys}) reduces to $D\ed V_0 = xy\ed
V_1,$ or in components:
\begin{align}\label{generic221}
  \left\{
    \begin{aligned}
      \frac{\partial V_0}{\partial x} = y\frac{\partial V_1}{\partial
        x}\\  
      \frac{\partial V_0}{\partial y} = x\frac{\partial V_1}{\partial
        y}.  
    \end{aligned}
  \right.
\end{align}    
This system has the general solution 
\begin{align*}
  V_0 = \frac{x\phi (y) - y\psi (x)}{x - y},\quad
  V_1 = \frac{\phi (y) - \psi (x)}{x - y}, 
\end{align*}
where $\phi$ and $\psi$ are arbitrary functions of one variable. 

\subsubsection{$(n,m,k) = (m,m,k)$}

We will now study the generic case of the more general situation when
the only restriction for the parameters $n,m,k$ is that $n = m.$ We
assume now that $q^1 = Z_0,\ldots ,q^n = Z_{m-1}$ are functionally
independent and constitute a complete set of coordinates. The
condition (\ref{chareqmatrix}), which guarantees existence of
multiplication for the system (\ref{sys}), now attains the simple form
$A_C = 0$ since $Z'$ in these coordinates becomes the identity matrix.
Hence, the system (\ref{sys}) admits $*-$multiplication if and only if
it, in the matrix notation with the coordinates $q^i = Z_{i-1},$ can
be written as
\begin{align}\label{generic_mmk}
  0 = \sum_{i=1}^k\left( C^iV' - V'C^i \right)A_i
\end{align}
where the tensors $A_1,\ldots , A_k$ are arbitrary. 

\subsubsection{$(n,m,k) = (m,m,1)$}

When $k=1,$ the equation (\ref{generic_mmk}) becomes a remarkably
simple equation $0 = (CV' - V'C)A_1$ in terms of the generic
coordinates. Thus, if we also assume that $A_1$ is non-singular, we
obtain in the generic $(m,m,1)-$case, the equation
\begin{align}\label{generic_mm1}
  CV' = V'C.
\end{align}
By calculating the residue of $A_\mu \ed V_\mu$ modulo $Z_\mu,$ we see
that the equation (\ref{generic_mm1}) can also be written as 
\begin{align}\label{generic_mm1b}
  -C\ed V_i + \ed V_{i-1} = q^{i+1}\ed V_{n-1}, \quad i = 0, 1, \ldots
  , n-1, \quad V_{-1} := 0
\end{align}  
However, the equation for which $i = 0$ in (\ref{generic_mm1b}) can be
discarded since it is a consequence of the other equations and of the
Cayley--Hamilton theorem for the companion matrix $C.$ This is
realized by adding to the equation in (\ref{generic_mm1b}) for which
$i = 0$ the equation for which $i = 1$ multiplied with $C,$ then
adding the equation for which $i = 2$ multiplied with $C^2,$ and so
on. In components, the equations (\ref{generic_mm1b}) can then be
written as
\begin{align*}
  \left\{
    \begin{aligned}
       -\pd{V_i}{q^2} + \pd{V_{i-1}}{q^1}  & =
       q^{i+1}\pd{V_{n-1}}{q^1}\\   
       -\pd{V_i}{q^3} + \pd{V_{i-1}}{q^2}  & =
       q^{i+1}\pd{V_{n-1}}{q^2}\\  
      & \;\; \vdots\\
       -\pd{V_i}{q^n} + \pd{V_{i-1}}{q^{n-1}}  & =
       q^{i+1}\pd{V_{n-1}}{q^{n-1}}\\   
       q^1\pd{V_i}{q^1} + \cdots + q^n\pd{V_i}{q^n} +
       \pd{V_{i-1}}{q^n}  & = q^{i+1}\pd{V_{n-1}}{q^n},  
    \end{aligned}
  \right. 
\end{align*} 
where $i = 1, \ldots , n-1.$ Thus, we see that the generic
$(m,m,1)-$case consists of $m(m-1)$ PDE's for $m$ dependent variables,
and is therefore an overdetermined system when $m>1$, that
nevertheless has non-trivial solutions, e.g., $V_i = Z_i = q^{i+1}.$
If we let $U = V*W,$ each entry $U_a$ of the $1-$column matrix $U$ is
a sum of terms $P_{ij}V_iW_j,$ where $P_{ij}$ is a polynomial of
degree at most $m-1$ in the coordinates $q^1,\ldots ,q^n.$ We note
also that for every solution $V$ of equation (\ref{generic_mm1}), each
term of the sum in equation (\ref{generic_mmk}) vanishes. Thus, every
solution $V$ in the generic $(m,m,1)-$case also solves the generic
$(m,m,k)-$case for arbitrary $k.$ The lowest order $*-$powers of the
trivial solution $\mu$ (or $[0,1,0,\ldots ,0]^T$ in matrix notation)
are given by:
\begin{align*}
  \mu^1_* & = \left[0,1,0,\ldots ,0\right]^T\\
  \mu^2_* & = \left[0,0,1,\ldots ,0\right]^T\\
  &\;\; \vdots\\
  \mu^{m-1}_* & = \left[0,0,\ldots ,1\right]^T\\
  \mu^m_* & = -\left[q^1,q^2,\ldots ,q^m\right]^T\\
  \mu^{m+1}_* & = \left[q^1q^m,-q^1+q^2q^m,\ldots ,-q^{m-1} + (q^m)^2\right]^T\\ 
  \mu^{m+2}_* & = \left[q^1(q^{m-1} - (q^m)^2), q^1q^m + q^2(q^{m-1} -
  (q^m)^2),\ldots ,\right.\\
  & \qquad \left. -q^{m-2} + q^m(2q^{m-1} - (q^m)^2) \right]^T
\end{align*}
\subsection{$m<n$}

Suppose that $Z_0, \ldots ,Z_{m-1}$ are functionally independent and
consider local coordinates $q^1,\ldots ,q^n$ such that $q^1 = Z_0,
\ldots , q^m = Z_{m-1},$ without specifying the other coordinates
$q^{m+1},\ldots ,q^n.$ Now (\ref{chareqmatrix}) is equivalent to $0 =
\sum_{i=0}^k\left[ C^i\;|\;0 \right] A_i,$ where $\left[
C^i\;|\;0\right]$ denotes a $m\times n-$matrix constructed by writing
$n-m$ zero-columns right to the matrix $C^i.$ This condition
determines the first $m$ rows of $A_0$ uniquely, while the last rows
as well as $A_1,A_2,\ldots ,A_k$ are arbitrary.\newline

One should note that when $m<n$, since $Z_i$ are the only functions
except $V_i$ and $W_i$ that appear in $V*W,$ this product will not
depend on certain coordinates unless $V$ or $W$ are themselves
functions depending on these coordinates. Thus, a product of two
trivial solutions will will never depend on these missing
coordinates.\newline

\subsubsection{$(n,m,k) = (n,2,1)$}

We will now choose the lowest possible values for $m$ and $k$ and let
$n$ be arbitrary. We note that the $(3,2,1)-$case, which is the
simplest possible case for which $m<n,$ is included. In generic
coordinates, the relation (\ref{chareqmatrix}) is equivalent to the
following relation between the components of the tensors $A_0$ and
$A_1:$
\begin{align*}
  \left[
  \begin{array}{ccc}
    (A_0)^1_1 & \cdots & (A_0)^1_n\\
    (A_0)^2_1 & \cdots & (A_0)^2_n
  \end{array}
  \right] = C
  \left[
  \begin{array}{ccc}
    (A_1)^1_1 & \cdots & (A_1)^1_n\\
    (A_1)^2_1 & \cdots & (A_1)^2_n
  \end{array}
  \right],
\end{align*}
where $C$ denotes the companion matrix of $Z_\mu.$ We see that, even
though we restrict our attention to the case $A_1=I,$ $n(n-2)$
components of $A_0$ can still be chosen arbitrary. However, for some
choices of these components, the system (\ref{sys}) will not depend on
some of the coordinates and can therefore be reduced to a lower
dimensional problem with a smaller number of independent variables. We
illustrate this phenomenon in the case $(3,2,1).$\newline

\subsubsection{$(n,m,k) = (3,2,1)$} 

According to the discussion above, we can assume that
\begin{align*}
  A_0 = \matrist{0}{x}{0}{-1}{y}{0}{a}{b}{c},
\end{align*}   
in generic coordinates $(x,y,z)$ where $a,b,c$ are arbitrary
functions. By analyzing the corresponding system (\ref{sys}), one can
see that either all solutions $V_0,$ $V_1$ are constant with respect
to the variable $z$ and the system then reduces to the generic
$(2,2,1)-$case, or otherwise we must have $y^2 > 4x,$ $b = a(c-y),$
and $c = y \pm (1/2)\sqrt{y^2-4x},$ and the system can then be written
as $A_0\ed V_0 = x\ed V_1,$ or in components
\begin{align*}
  \left\{
    \begin{aligned}
      & -\pd{V_0}{y} + a(x,y,z)\pd{V_0}{z} = x\pd{V_1}{x}\\    
      & x\pd{V_0}{x} + y\pd{V_0}{y} +
      a(x,y,z)\left(\frac{-y\pm\sqrt{y^2-4x}}{2}\right)\pd{V_0}{z} =
      x\pd{V_1}{y}\\  
      & \left(\frac{y\pm\sqrt{y^2-4x}}{2}\right)\pd{V_0}{z} =
      x\pd{V_1}{z},  
    \end{aligned}
  \right.
\end{align*}     
where $a = a(x,y,z)$ is an arbitrary function. Thus, we see that the
generic $(3,2,1)-$case involves an arbitrary function, which was not
possible for the $(m,m,1)$-case. We note also that all solutions $V$
of the generic $(2,2,1)$-case also solve the generic $(3,2,1)-$case.
Since the multiplication coincides with the multiplication in the
$(2,2,1)-$case, we see that unless $V$ or $W$ depend on $z,$ the
product $V * W$ will not depend on $z.$ Especially, the solutions
obtained by taking powers of the trivial solution $\mu$ are again
given by (\ref{tp}).\newline

\subsubsection{$(n,m,k) = (n,3,2)$}

When we consider higher values of the number $k,$ the corresponding
systems become harder to analyze. One reason is that when we increase
$k$ by one, we add a new tensor $A_k$ which means that we add
$n^2$ new components. Another reason is that with higher values of
$k,$ we get higher order polynomials in $Z_i.$ Already for $k =
2,$ such systems become quite hard to handle. In the general
$(n,3,2)-$case for example, the system (\ref{sys}) can be written as
\begin{align*}
  \left\{
    \begin{aligned}
      0  = &\; A_0\ed V_0 - Z_0A_1\ed V_2 - Z_0A_2\ed V_1 +
      Z_0Z_2A_2\ed V_2\\  
      0  = &\; A_0\ed V_1 + A_1\ed V_0 - Z_1A_1\ed V_2 - Z_1A_2\ed V_1
      + \left(Z_1Z_2 - Z_0\right)A_2\ed V_2\\ 
      0  = &\;  A_0\ed V_2 + A_1\ed V_1 + A_2\ed V_0 - Z_2A_1\ed V_2 -
      Z_2A_2\ed V_1 + \left( Z_2^2 - Z_1 \right)A_2\ed V_2.   
    \end{aligned}
  \right.
\end{align*}    
Even if we assume that $A_2=I$ and that $q^{i+1} = Z_i$ are
functionally independent, we still have $2n^2$ arbitrary functions in
the picture (the components of $A_0$ and $A_1$). When the condition
(\ref{chareqmatrix}), which in generic coordinates is equivalent to $0
= \sum_{i=1}^k\left[ C^i\;|\;0\right] A_i,$ is satisfied we still have
$n(2n-m)$ arbitrary functions.

\subsection{$m > n$}

When $m > n$, the generic case becomes more complicated than for $m
\le n.$ Consider the generic case when $q^1 = Z_0,\ldots ,q^n =
Z_{n-1}$ are functionally independent and take $q^1, \ldots ,q^n$ as
local coordinates. In these coordinates we have
\begin{align*}
  Z'=\left[
  \begin{array}{cccc}
    1 & 0 & \cdots & 0 \\
    0 & 1 & \cdots & 0 \\
    \vdots & & \ddots & \vdots \\
    0 & 0 & \cdots & 1 \\
    \partial_1 Z_n & \partial_2 Z_n & \cdots & \partial_n Z_n \\ 
    \vdots & \vdots & & \vdots \\
    \partial_1 Z_{m-1} & \partial_2 Z_{m-1} & \cdots & \partial_n
    Z_{m-1} \\ 
  \end{array}
  \right].
\end{align*} 
Thus, the relation (\ref{chareqmatrix}), which in the generic cases
for $m \le n$ became a set of algebraic equations for the components
of $A_i,$ becomes now a complicated differential relation between the
components of $A_i$ and the functions $Z_n,\ldots ,Z_{m-1}.$

\subsubsection{$(n,m,k) = (2,3,1)$}

We consider the simplest case for which $m > n,$ i.e., when $(n,m,k) =
(2,3,1).$ The system (\ref{sys}) can then be written as
\begin{align}\label{231sys}
  \left\{
  \begin{aligned}
    0 & = \; A_0\ed V_0 - Z_0A_1\ed V_2\\
    0 & = \;A_1\ed V_0 + A_0\ed V_1 - Z_1A_1\ed V_2\\ 
    0 & = \; A_1\ed V_1 + \left(A_0 - Z_2A_1\right)\ed V_2. 
  \end{aligned}
  \right.
\end{align}  
The condition $m > n,$ implies that the functions $Z_0,$ $Z_1,$ $Z_2$
are functionally dependent. For the sake of convenience, we assume in
the generic case that the functions $Z_0,$ $Z_1$ are functionally
independent and that $Z_2 = \phi (Z_0,Z_1)$ for some function $\phi.$
We also assume that $A_1 = I.$ Thus, in generic coordinates $x = Z_0,$
$y = Z_1,$ the condition (\ref{chareqmatrix}) is satisfied if and only
if $A := A_0$ is given by
\begin{align*}
  A = \matris{x\partial_x \phi}{x\partial_y
    \phi}{y\partial_x \phi - 1}{y\partial_y \phi}, 
\end{align*}  
where $\phi$ satisfies the non-linear partial differential
equations  
\begin{align}\label{231pde}
  \left\{
  \begin{aligned}
    0 = & \;x\left(\pd{\phi}{x}\right)^2 - \phi\pd{\phi}{x} +
    y\pd{\phi}{x}\pd{\phi}{y} - \pd{\phi}{y}\\ 
    0 = & \;1 + x\pd{\phi}{x}\pd{\phi}{y} + y\left( \pd{\phi}{y}
    \right)^2 - \phi\pd{\phi}{y}. 
  \end{aligned}
  \right.
\end{align}  
Solutions of the equations (\ref{231pde}) exist, for instance $\phi =
ay - a^2x + a^{-1}$ solves (\ref{231pde}) for any non-zero real
constant $a.$ In the generic case, the system (\ref{231sys}) is
equivalent to (note that the tensor $A$ is non-singular when $x\ne 0$) 
\begin{align*}
  \left\{
  \begin{aligned}
    0 = & \;A\ed V_0 - x\ed V_2\\
    0 = & \;\left( xI - yA + \phi A^2 - A^3 \right)\ed V_2\\
    0 = & \;\ed V_1 + \left(A - \phi I \right) \ed V_2.  
  \end{aligned}
  \right.
\end{align*}  
The Cayley--Hamilton theorem together with the assumption that $\phi$
satisfies the equations (\ref{231pde}) implies that the expression $xI
- yA + \phi A^2 - A^3$ vanishes. Hence, the system (\ref{231sys})
reduces to 
\begin{align*}
  \left\{
  \begin{aligned}
    A\ed V_0 & = \;x\ed V_2\\
    \ed V_1 & =  \;\left(\phi I - A \right) \ed V_2.  
  \end{aligned}
  \right.
\end{align*}  
In the case when $\phi = ay - a^2x + a^{-1},$ the system
(\ref{231sys}) can explicitly be written as the overdetermined system   
\begin{align*}
  \left\{
  \begin{aligned}
   0 & = \;a^2x\pd{V_0}{x} + (ay^2 + 1)\pd{V_0}{y} + x\pd{V_2}{x} \\ 
   0 & = \;-ax\pd{V_0}{x} - ay\pd{V_0}{y} + x\pd{V_2}{y}\\ 
   0 & = \;-\pd{V_1}{x} + (ay + a^{-1})\pd{V_2}{x} +
   (a^2y+1)\pd{V_2}{y}\\   
   0 & = \;-\pd{V_1}{y} - ax\pd{V_2}{x} + (a^{-1} - a^2x)\pd{V_2}{y}. 
  \end{aligned}
  \right.
\end{align*}  
For this system, the $*-$product of trivial solutions gives in general
non-trivial solutions. For example, when $V = [0, 1, 0]^T$ and $W =
[0, 0, 1]\;^T$ (or simply $V_\mu = \mu$ and $W_\mu = \mu_*^2$ using
the $\mu-$notation), their $*-$product becomes
\begin{align*}
  V*W = \vektort{-x}{-y}{a^2x-ay-a^{-1}}.
\end{align*} 

\subsection{Summary}

The results about the possible structures of the system of linear
PDE's (\ref{sys}), and of the corresponding relation (\ref{chareq})
(that characterizes the existence of $*-$multi-\linebreak plication),
can for non-singular $A_k$ be summarized as follows:
\begin{enumerate}
  \item $m=n.$ Suppose that $Z_0, \ldots ,Z_{m-1}$ are
    functionally independent. In the generic coordinates $q^i =
    Z_{i-1},$ the relation (\ref{chareq}) is satisfied if and only
    if $A_C = A_0C^0 + \cdots A_kC^k = 0,$ and the system
    (\ref{matrix_sys}) can be written as 
    \begin{align*} 
      0 = \sum_{i=1}^k\left( C^iV' - V'C^i \right)A_i, 
    \end{align*}  
    where the tensors $A_1, \ldots ,A_k$ are arbitrary. In the
    generic case when $k=1,$ the system (\ref{matrix_sys}) takes the
    remarkably simple form $CV' = V'C,$ for which every solution $V$
    is also a solution in each generic $(n,m=n,k)-$case with arbitrary
    $k>0$. 
  \item $m<n.$ Suppose that $Z_0, \ldots ,Z_{m-1}$ are
    functionally independent. In the generic coordinates $q^1 =
    Z_0, \ldots , q^{m} = Z_{m-1}, q^{m+1}, \ldots q^n$
    ($q^{m+1}, \ldots ,q^n$ not specified), the relation
    (\ref{chareq}) is satisfied if and only if $0 = \sum_{i=0}^k\left[
      C^i\;|\;0\right] A_i.$  It is characteristic for the case $m<n$
    that $*-$products of trivial solutions will remain constant in
    some variables.  
  \item $m>n.$ This is the hardest case to analyze since the relation
    (\ref{chareq}) in the generic coordinates $q^1 = Z_0, \ldots ,
    q^{n} = Z_{n-1}$ becomes a differential relation for the
    functions $Z_{n+1}, \ldots , Z_{m-1},$ while it leads to
    algebraic equations in the previous cases.
\end{enumerate}

\section{Power series}\label{sec_ps}

As mentioned above, we can take power series of the trivial solution
$\mu$ to build up more complicated solutions of (\ref{sys}). In the
matrix notation, such power series have the form
\begin{align*}
  P = \sum_{r=0}^\infty a_rC^re_1,
\end{align*} 
where $a_r$ are real constants. We will now investigate when these
power series define new solutions of the system (\ref{sys}), i.e. when
they are convergent and the summation and derivation commutes so that
they define genuine solutions of the first order systems of PDE's
(\ref{sys}).\newline 

The companion matrix $C$ can be factorized as $C = TJT^{-1},$ where
the matrix $J$ has the Jordan canonical form, i.e.
$J=\mathrm{diag}(J_1,\ldots ,J_s),$ where
\begin{align}\label{jordanblock}
  J_s = \left[ 
    \begin{array}{cccc}
      \lambda_s & 1 & &\\
      & \lambda_s & \ddots &\\
      & & \ddots & 1\\
      & & & \lambda_s
    \end{array} \right]
\end{align} 
and $\lambda_s$ is an eigenvalue of $C$ (note that the eigenvalues of
$C$ coincide with the roots of the polynomial $Z_\mu$). Thus, the
partial sums of the power series can be written as 
\begin{align}
  P_N = \sum_{r=0}^N a_rC^re_1 = T\left( \sum_{r=0}^N a_rJ^r \right)
  T^{-1} e_1. \label{pjordan}
\end{align}

\begin{TL}\label{jpower}
  Let $J_s$ be the Jordan block defined by $($\ref{jordanblock}$)$. The
  entries of powers of $J_s$ can are given by 
  \begin{align*}
    \left( J^r_s \right)_{ij} = {r \choose j-i}\lambda^{r+i-j}_s. 
  \end{align*}
\end{TL}

\begin{proof}

  The proof follows immediately by induction over $r.$

\end{proof}
From the lemma, it is clear that all non-zero elements of the matrix
$\sum_{r=0}^N a_rJ^r$ have the form $\sum_{r=0}^N a_r{r \choose
b}\lambda_s^{r-b},$ for eigenvalues $\lambda_s$ of $C.$ All power
series \linebreak $\sum_{r=0}^\infty a_r{r \choose b}\lambda_s^{r-b}$ have the
same radius of convergence as $\sum_{r=0}^\infty a_r t^r.$ Hence, we
conclude that at every point of $Q$ where all eigenvalues of $C$
belongs to the open set $]-R,R[$ (where $1/R = \limsup |a_n|^{1/n}$),
the power series $\sum_{r=0}^\infty a_rJ^r,$ and therefore also
$\sum_{r=0}^\infty a_rC^re_1,$ will be convergent. The question that
now remains, is what is required for a convergent power series
$\sum_{r=0}^\infty a_rC^re_1$ to be a solution to (\ref{sys}). The
following theorem gives a sufficient condition:

\begin{T}\label{pst}
  Consider a power series $P = \sum_{r=0}^\infty a_rC^re_1.$
  Let $D\subset Q$ be a domain in which each eigenvalue $\lambda_s$ of
  the matrix $C$ is real and $\lambda_s\in [-R+\epsilon,R-\epsilon ]$
  for some $\epsilon > 0,$ where $1/R = \limsup|a_r|^{1/r}.$ Assume
  also that the geometrical multiplicity of each eigenvalue is
  constant in $D.$ Then the power series $P$ defines a solution of
  $($\ref{sys}$)$ in $D.$
\end{T}

\begin{Remark}
  The assumption about the geometrical multiplicity guarantees that
  the structure of the Jordan canonical form of $C$ is preserved in
  $D$ (note that the case with simple eigenvalues is included). There
  is no indication that this assumption, or the assumption about
  eigenvalues being real, are necessary for the conclusion in theorem
  \ref{pst}, but without them the proof would be technically more
  complicated. 
\end{Remark} 

\begin{proof}
  $P$ is a solution of (\ref{sys}) if and only if 
  \begin{align*}
    0 = \sum_{i=0}^k C^iP'A_i = \sum_{i=0}^k C^i\left(
      \sum_{r=0}^\infty a_rC^re_1 \right)'A_i.  
  \end{align*}
  Thus, since $C^re_1$ is a solution for any $r,$ we see that it is
  enough to prove that $P' = \lim_N P_N',$ or in components,
  $\partial_i P^j = \lim_N \partial_i P_N^j,$ $i,j=1,\ldots ,n,$ where
  $P^j$ and $P^j_N$ denote the $j'$th element of $P$ and $P_N,$
  respectively. Since $P_N$ converges to $P$ in every point, we only
  have to prove that $\partial_i P^j_N$ converges uniformly in $D.$
  \newline 

  According to (\ref{pjordan}) and lemma \ref{jpower}, each $P_N^j$
  consists of a finite sum of terms of the form $f(q)\sum_r^N
  b_r\lambda^r_s,$ where $\lambda_s$ is an eigenvalue of $C,$ $f$ is a
  real-valued function, and the power series $\sum_r^\infty
  b_r\lambda^r_s$ has $R$ as its radius of convergence. Thus,
  $\partial_iP_N^j$ consists of a finite sum of terms of the form
  $(\partial_i f(q))\sum_r^N b_r\lambda^r_s + f(q)(\partial_i\lambda_s
  )\sum_r^N rb_r\lambda^{r-1}_s,$ which converge uniformly in $D.$ 
\end{proof}

According to theorem \ref{pst}, if the companion matrix $C$ has
simple eigenvalues and $\sum_r a_rt^r$ is a power series with infinite
radius of convergence, then $\sum_r a_rC^re_1$ is a globally defined
power series solution of (\ref{sys}). Thus, for example, we can
construct the power series solutions
\begin{align*} 
  \mathrm{exp}_* C & := \left( \sum_{r=0}^\infty \frac{1}{r!}C^r
  \right)e_1\\   
  \mathrm{sin}_* C & := \left( \sum_{r=1}^\infty
    (-1)^{r-1}\frac{1}{(2r-1)!}C^{2r-1} \right)e_1\\ 
  \mathrm{cos}_* C & := \left( \sum_{r=0}^\infty
    (-1)^r\frac{1}{(2r)!}C^{2r} \right)e_1. 
\end{align*}

We end the discussion about power series with some examples.

\begin{Example}

  To illustrate the mechanism of generating power series solutions,
  with respect to the $*-$multiplication, we return to the system
  $($\ref{sys221}$)$ in the $(2,2,1)-$\linebreak case. We assume that
  $A_\mu$ and $Z_\mu$ satisfy the relation $($\ref{chareq}$)$ so that
  the system has $*-$multiplication, and we also assume that the roots
  $\lambda_1,\; \lambda_2,$ of $Z_\mu$ are simple. The companion
  matrix can then be diagonalized as 
  \begin{align*}
    C = TDT^{-1} = \matris{\lambda_2}{\lambda_1}{-1}{-1}
    \matris{\lambda_1}{0}{0}{\lambda_2}
    \matris{\lambda_2}{\lambda_1}{-1}{-1}^{-1}. 
  \end{align*}
  A power series $\sum_r a_rC^re_1$ defines a solution in any domain
  $D\subset Q$ in which $|\lambda_1|,\; |\lambda_2| <
  1/\limsup|a_r|^{1/r} - \epsilon,$ and can be written as
  \begin{align*}
    \sum_r a_rC^re_1 & = \;T\left( \sum_r a_rD^r \right) T^{-1}e_1\\
    & = \;\frac{1}{\lambda_1 - \lambda_2}
    \matris{\lambda_2}{\lambda_1}{-1}{-1} \matris{\sum
      a_r\lambda_1^r}{0}{0}{\sum a_r\lambda_2^r}
    \matris{-1}{-\lambda_1}{1}{\lambda_2} \vektor{1}{0}\\
    & = \;\frac{1}{\lambda_1 - \lambda_2} \vektor{\lambda_1\sum
      a_r\lambda_2^r - \lambda_2\sum a_r\lambda_1^r}{\sum
      a_r\lambda_1^r - \sum a_r\lambda_2^r}.
  \end{align*}
  Thus, for example we have
  \begin{align*}
    \mathrm{exp}_* C & = \frac{1}{\lambda_1 - \lambda_2}
    \vektor{\lambda_1 \exp{(\lambda_2)} - \lambda_2
      \exp{(\lambda_1)}}{\exp{(\lambda_1)} - \exp{(\lambda_2)}}.  
  \end{align*}  
  In the case when the eigenvalues of $C$ are constant (which is the
  case for the Cauchy--Riemann equations), these power series will
  only provide trivial solutions. Hence, in order to obtain
  interesting solutions for this case, one has to build power series solutions from
  a non-trivial solution (for instance $(x,y)$ in the Cauchy--Riemann
  case). The other extreme case is the generic situation when the
  eigenvalues are functionally independent and $A_1$ is
  non-singular. Then, as mentioned earlier, the general solution of
  $($\ref{generic221}$)$ can be written as    
  \begin{align}\label{221sol}
    V = \frac{1}{\lambda_1 - \lambda_2} \vektor{ \lambda_1\phi
      (\lambda_2) - \lambda_2\psi (\lambda_1)}{\psi (\lambda_1) - \phi
      (\lambda_2)}.   
  \end{align}
  It is a remarkable property that every analytic solution in the
  generic case can be expressed by power series of the trivial solution
  $(0,1)$. Namely, if $V$ is a solution of the form $($\ref{221sol}$)$
  where $\phi(t) = \sum a_rt^r$ and $\psi(t)=\sum b_rt^r$ are analytic,
  then $V$ can be expressed in terms of the two simple solutions
  $(\lambda_1 - \lambda_2)^{-1}(\lambda_1, -1)$ and $(\lambda_1 -
  \lambda_2)^{-1}(-\lambda_2, 1),$ and of power series of the trivial
  solution $(0,1)$:
  \begin{align*}
    V = \frac{1}{\lambda_1 - \lambda_2}\vektor{\lambda_1}{-1} *
    \left( \sum a_rC^re_1 \right) + \frac{1}{\lambda_1 -
      \lambda_2}\vektor{-\lambda_2}{1} * \left( \sum b_rC^re_1 \right)
  \end{align*}
  Hence, a significant part of the solution set of
  $($\ref{generic221}$)$ can be expressed by power series in trivial
  solutions. 

  \begin{Remark}
    The solutions $(\lambda_1 - \lambda_2)^{-1}(\lambda_1, -1)$ and
    $(\lambda_1 - \lambda_2)^{-1}(-\lambda_2, 1)$ have some
    remarkable properties. They are idempotent and their sum is the
    identity, i.e.,
    \begin{align*}
      \begin{aligned}
        &\left(\frac{\lambda_1}{\lambda_1-\lambda_2},
          \frac{-1}{\lambda_1-\lambda_2}\right)_*^2 =
        \left(\frac{\lambda_1}{\lambda_1-\lambda_2},
          \frac{-1}{\lambda_1-\lambda_2}\right),  \\  
        &\left(\frac{-\lambda_2}{\lambda_1-\lambda_2},
          \frac{1}{\lambda_1-\lambda_2}\right)_*^2 =
        \left(\frac{-\lambda_2}{\lambda_1-\lambda_2},
          \frac{1}{\lambda_1-\lambda_2}\right), \\ 
        &\left(\frac{\lambda_1}{\lambda_1-\lambda_2},
          \frac{-1}{\lambda_1-\lambda_2}\right) +
        \left(\frac{-\lambda_2}{\lambda_1-\lambda_2},
          \frac{1}{\lambda_1-\lambda_2}\right) = (1,0).
      \end{aligned}  
    \end{align*}
  \end{Remark}
\end{Example}

\begin{Example}
  Let us instead consider the case where $m = 2$ and the eigenvalues
  of $C$ coincide everywhere, $\lambda_1 = \lambda_2 =: \lambda$, i.e.,
  the functions $Z_0$ and $Z_1$ are functionally dependent and related
  as $Z_1^2 = 4Z_0.$ Then we can factorize $C$ as
  \begin{align*}
    C = TJT^{-1} = \matris{1}{-\lambda}{0}{1}
    \matris{\lambda}{0}{1}{\lambda} \matris{1}{\lambda}{0}{1},  
  \end{align*}
  and a general power series solution can be written as
  \begin{align*}
    \sum_r a_rC^re_1 = T \left( \sum_r a_rJ^r \right) T^{-1}e_1 =
    \vektor{\sum_r (1-r)a_r\lambda^r }{ \sum_r ra_r\lambda^{r-1} }. 
  \end{align*}
  For this example, the exponential power series produces the solution
  $\mathrm{exp}_* C = \exp{(\lambda)(1 - \lambda, 1)}.$  
\end{Example}

\begin{Example}\label{jodeitex}
  In \cite{jodeit-1990}, a matrix equation of the form
  \begin{align}\label{fmg}
    \nabla f = M \nabla g,
  \end{align}
  where $M$ is a constant matrix, is studied in an open convex domain
  of a real or complex vector space. Here $f(\mathbf{x})$ and
  $g(\mathbf{x})$ are scalar-valued functions defined on the vector
  space, and $\nabla f = [\partial_{x_1}f,\ldots ,\partial_{x_n}f]^T$,
  where $\mathbf{x} = [x_1,\ldots ,x_n]^T$ is a coordinate vector with
  respect some basis. We give below an example of an equation of the
  form $($\ref{fmg}$)$ that admits a $*-$multiplication, and we show that
  solutions can be represented by power series of simple solutions.
  \newline

  Let $V$ be a real vector space of dimension four, and assume that
  the matrix $M$ in $($\ref{fmg}$)$ has eigenvalues $\lambda$, $\lambda$,
  $\alpha \pm i\beta$, where $\lambda$, $\alpha$, and $\beta$ are real
  constants. By performing a linear change of variables $\mathbf{x}
  \rightarrow A\mathbf{x}$, where $A$ is a constant matrix, the
  equation $($\ref{fmg}$)$ transforms to $\nabla f = A^{-T}MA^T \nabla g$.
  Thus, we can assume that $M$ is in canonical real normal form. We
  assume that the eigenvalue $\lambda$ has geometrical multiplicity
  one, so that $M$ is given by
  \begin{align}\label{cstm}
    M = \left[
      \begin{array}{cccc}
        \lambda & 1 & 0 & 0 \\
        0 & \lambda & 0 & 0 \\
        0 & 0 & \alpha & \beta \\
        0 & 0 & -\beta & \alpha
      \end{array}
    \right].
  \end{align}

  According to \cite{jodeit-1990}, the general analytic solution of
  $($\ref{fmg}$)$ can be decomposed as $f = f_1 + f_2,$ $g = g_1 + g_2$,
  where $(f_1,\;g_1)$ and $(f_2,\;g_2)$ are solutions of 
  \begin{align}\label{subsys}
    \nabla f_1 = \matris{\lambda}{1}{0}{\lambda}\nabla g_1,\quad and
    \quad \nabla f_2 =
    \matris{\alpha}{\beta}{-\beta}{\alpha} \nabla g_2,  
  \end{align}
  respectively. By changing to new dependent variables
  \begin{align*}
    \tilde{f}_1 = f_1 - \lambda g_1, \quad \tilde{g}_1 = g_1, \quad
    \tilde{f}_2 = f_2 - \alpha g_2, \quad \tilde{g}_2 = \beta g_2, 
  \end{align*}
  it becomes obvious that there is no restriction to assume that
  $\lambda = \alpha = 0$ and $\beta = 1$. We note that the equation
  for $f_2$ and $g_2$ then reduces to the Cauchy--Riemann equations.
  For this choice of $M$, it is trivial to obtain the general analytic
  solution of $($\ref{fmg}$)$:     
  \begin{align*}
    \left\{
      \begin{aligned}
        f_1 & = x_2\phi'(x_1) + \psi(x_1) \\
        g_1 & = \psi(x_1) + c
      \end{aligned}
    \right.\quad
    \left\{
      \begin{aligned}
        f_2 & = \mathrm{Re}\left( F(x_3 + \mathrm{i}x_4) \right)\\
        g_2 & = \mathrm{Im}\left( F(x_3 + \mathrm{i}x_4) \right),
      \end{aligned}
    \right.
  \end{align*}
  where $\psi$ and $\phi$ are arbitrary analytic functions of one real
  variable, and $F$ is an arbitrary holomorphic function of one
  complex variable. We have already seen that the general solution of
  the Cauchy--Riemann equations can be represented as a power series
  of the simple solution $x_3 + \mathrm{i}x_4$ with respect to a
  $*-multiplication$ (coinciding with multiplication of holomorphic
  functions). Also the first system in $($\ref{subsys}$)$ admits a
  $*-$multiplication since it can be written as a system of the form
  $($\ref{sys}$)$ where
  \begin{align*}
    A_\mu = \matris{0}{0}{-1}{0} + \mu I,\quad V_\mu = f_1 + \mu
    g_1,\quad Z_\mu = \mu^2,  
  \end{align*}
  and the relation $($\ref{chareq}$)$ is trivially satisfied since
  $Z_\mu$ is constant. It is remarkable that also for this system,
  the general analytic solution can written as a power series of
  simple solutions. Namely, let the functions $\phi$ and $\psi$, from 
  the general solution, be analytic with power series representations
  $\psi(s) = \sum_r a_rs^r$ and $\phi(s) = \sum_r b_rs^r$,
  respectively. Then, the solution $(f_1,g_1)$ is given by
  \begin{align*}
    f_1 + \mu g_1 = \sum_r (a_r + \mu b_r)*(x_1 + \mu x_2)^r_*.
  \end{align*}
  Thus, when $M$ is given by $($\ref{cstm}$)$, we have shown that the
  general analytic solution of $($\ref{fmg}$)$ can be obtained by taking power
  series, with respect to $*-$multiplication, of simple solutions of
  the subsystems $($\ref{subsys}$)$. Also for general systems of type
  $($\ref{fmg}$)$, solutions can be built up by taking  $*-$power
  series of simple non-trivial solutions.
\end{Example}

\section{How to find systems with multiplication}\label{sec_find}

As we have seen, the problem of finding systems of the form
(\ref{sys}) that allow a $*-$multiplication is equivalent to finding a
$(1,1)-$tensor $A_\mu$ and a function $Z_\mu$ such that $A_\mu\ed
Z_\mu\equiv 0.$ In this general form, the problem is hard to handle,
since in coordinates we may need to solve a system of
complicated non-linear PDE's.\newline

One way to construct systems of PDE's allowing multiplication, is to
choose the function $Z_\mu$ first and treat it as fixed. The equation
(\ref{chareq}), or equivalently (\ref{chareqmatrix}), becomes then a
system of linear algebraic equations for the components of $A_i,$
which is easier to solve. It constitutes a system of $mn$ equations
for $n^2(k+1)$ unknown functions. The number $n(k+1)-m$ decides about
how large family of solutions that can be found for each $Z_\mu.$ When
$n(k+1)-m<0,$ the number of unknown functions is less than the number
of equations and we may not expect to find solutions for every choice
of $Z.$ On the other hand, for large values of $n(k+1)-m$ we get many
families of solutions whenever the system is consistent. We illustrate
this process of finding systems with a $*-$multiplication by an
example. \newline

\begin{Example}

  Let $Q$ be a manifold of dimension two with local coordinates $x,y$
  and suppose that $Z_\mu = x + y\mu + xy\mu^2 + \mu^3$ and $k = 1.$
  When $Z_\mu$ is given, we can find all tensors $A_\mu = A_0 +
  \mu A_1$ such that the corresponding system $($\ref{sys}$)$ is
  equipped with a $*-$multiplication, i.e, such that $($\ref{chareq}$)$ is
  satisfied.  
  \begin{align*}
    0 & = Z'A_0 + CZ'A_1\\
    & = \left[\begin{array}{ccc}
        1 & 0\\
        0 & 1\\
        y & x
      \end{array}\right]
    \matris{(A_0)^1_1}{(A_0)^1_2}{(A_0)^2_1}{(A_0)^2_2} +\\
    & \quad \matrist{0}{0}{-x}{1}{0}{-y}{0}{1}{-xy} \left[\begin{array}{ccc}  
        1 & 0\\
        0 & 1\\
        y & x
      \end{array}\right]
    \matris{(A_1)^1_1}{(A_1)^1_2}{(A_1)^2_1}{(A_1)^2_2}.
  \end{align*}
  The solution of this system of linear equations is given by
  \begin{align*}
    A_0 = \matris{xy}{x^2}{y^2-1}{xy}A_1,\quad A_1 =
    \matris{(x^2y+1)f}{(x^2y+1)g}{x(1-y^2)f}{x(1-y^2)g} 
  \end{align*}
  where $f = f(x,y)$ and $g = g(x,y)$ are two arbitrary functions. The
  system $($\ref{sys}$)$ takes a form that is independent of the choice
  of $f$ and $g$:
  \begin{align}\label{findex}
    \left\{
      \begin{aligned}
        0 & = x(y+x^2)\pd{V_0}{x} + (y^2 - 1)\pd{V_0}{y} -
        x(1+x^2y)\pd{V_2}{x} + x^2(y^2 - 1)\pd{V_2}{y}\\
        0 & = (1+x^2y)\pd{V_0}{x} + x(1 - y^2)\pd{V_0}{y} + x(y
        + x^2)\pd{V_1}{x} + (y^2 - 1)\pd{V_1}{y} -\\
        & \qquad y(1 + x^2y)\pd{V_2}{x} + xy^2(y^2 -1)\pd{V_2}{y}\\
        0 & = (1 + x^2y)\pd{V_1}{x} + x(1 - y^2)\pd{V_1}{y} +
        x^3(1 - y^2)\pd{V_2}{x} + \\ & \qquad(y^2 - 1)(x^2y + 1)\pd{V_2}{y}
      \end{aligned}
    \right.
  \end{align}
  Although there are two arbitrary functions $f$ and $g,$ every choice
  give rise to the same system. Thus, for this particular choice of
  $Z_\mu$ and parameters $(n,m,k),$ $($\ref{findex}$)$ is the only
  system with $*-$multiplication. Since $m = 3$, the corresponding
  multiplication formula is given by $($\ref{mult3}$)$. 
\end{Example}

In \cite{jonasson-2006}, the problem of finding systems with
multiplication on Riemannian manifolds was studied for the special case
when $A_\mu = X +\mu I$ and $Z_\mu = \det{A_\mu}.$ With those
restrictions, the problem of finding a system equipped with a
multiplication reduces to finding a tensor $X$ that satisfies the
equation (\ref{xmueq2}). The following families of solutions were found in
\cite{jonasson-2006}: 

\begin{enumerate}
  \item $X=\jt^{-1}J$ where $J$ and $\jt$ are arbitrary special
    conformal Killing tensors. In this case, the multiplication of
    cofactor pair systems \cite{lundmark-2001} is reconstructed.

  \item Every $X$ with a vanishing Nijenhuis torsion $N_X = 0$. This
    follows from the remarkable relation $2\left(X\ed(\det{X}) -
      \det{X}\ed(\tr X)\right)_i = (N_X)^k_{ij}C^j_k,$ where $C = \cof
    X,$ and that $X$ and $X + \mu I$ share the same torsion. This
    result holds also when no metric is specified on the manifold $Q.$ 

  \item If $X$ is a non-singular solution, then $X^{-1}$ is a solution
    as well.

  \item In \cite{jonasson-2006}, a method for constructing solutions
    $X$ consisting of smaller blocks that satisfy (\ref{xmueq2}) is
    presented. A similar result is valid for the more general equation
    (\ref{chareq}):

    \begin{T}
      Suppose that $A_\mu$ and $Z_\mu$ satisfies the relation
      $($\ref{chareq}$)$ and let $\at_\mu$ be a $(1,1)-$tensor, and
      $\zt_\mu$ a function on a different manifold $\tilde{Q}$,
      satisfying the relation $\at_\mu\ed\zt_\mu \equiv 0 \quad
      (mod\;\zt_\mu).$ Then the tensor $A_\mu \oplus \at_\mu$ is a
      $(1,1)-$tensor on the manifold $Q\times \tilde{Q}$ and it
      satisfies the relation 
      \begin{align*}
        (A_\mu \oplus \at_\mu) \ed (Z_\mu\zt_\mu) \equiv 0 \quad
        (mod\;Z_\mu\zt_\mu). 
      \end{align*}
    \end{T}
    
    \begin{proof}
      \begin{align*}
        & (A_\mu \oplus \at_\mu) \ed (Z_\mu\zt_\mu) = (A_\mu \oplus
        \at_\mu) \left( \ed_Q (Z_\mu\zt_\mu) + \ed_{\tilde{Q}}
        (Z_\mu\zt_\mu) \right) \\
        &= A_\mu\ed_Q (Z_\mu\zt_\mu) + \at_\mu\ed_{\tilde{Q}} (Z_\mu\zt_\mu) =
        \zt_\mu A_\mu\ed_Q Z_\mu + Z_\mu\at_\mu\ed_{\tilde{Q}} \zt_\mu\\ 
        &= Z_\mu\zt_\mu\alpha \equiv 0 \quad (mod\;Z_\mu\zt_\mu),
      \end{align*}      
      where $\alpha$ is a certain $1-$form, and $\ed_Q$ and
      $\ed_{\tilde{Q}}$ denote the exterior differential operator on
      $Q$ and $\tilde{Q}$, respectively. 
    \end{proof}

  \item If the non-singular tensor $X$ is diagonal in the coordinates
    $q^1, \ldots ,q^n,$ then $X$ satisfies the equation (\ref{xmueq2})
    if and only if $X = \mathrm{diag}(X_1, \ldots , X_s),$ where
    $X_a = \mathrm{diag}(\phi_a, \ldots , \phi_a)$ is a square
    diagonal matrix of size $n_a$ and $\phi_a = \phi_a(q^{n_1 + \cdots
      + n_{a-1} + 1}, \ldots , q^{n_1 + \cdots + n_a})$ is an
    arbitrary (sufficiently regular) function, depending only on the
    specified coordinates.  

\end{enumerate}  

As we have seen, there are several ways to construct system of partial
differential equations with multiplication of solutions. Nevertheless,
for any given system of PDE's, it is in general hard to settle whether
it admits a $*-$multiplication.

\section{Conclusions}\label{sec_conc}

The $*-$multiplication constitutes a powerful method for generating,
in a pure algebraic way, new solutions from known solutions of certain
linear systems of PDE's. Especially by taking trivial solutions, we can
construct large families of non-trivial solutions of systems for which
non-trivial solutions are hard to obtain by other methods. By
generalizing the ideas from \cite{jonasson-2006}, we have
significantly extended the class of systems of PDE's admitting
$*-$multiplication. By identifying which elements of the construction
of systems of PDE's with multiplication that are relevant, we have
obtained much better understanding of the nature of the
$*-$multiplication. Our insight into the mechanism of the
multiplication has been obtained due to the matrix formulation of the
problem and of encoding the Euclidean algorithm as a matrix polynomial
(section \ref{matrixnot}), and due to the effective construction of
power series solutions (section \ref{sec_ps}). \newline

There are still many questions regarding the $*-$multiplication which
it may be worth to study. Some examples of such questions are:

\begin{enumerate}
  \item Which solutions can be represented as power series of
    trivial solutions, and when do they constitute all
    solutions? What can we say about power series of non-trivial but
    simple solutions (compare with the Cauchy--Riemann equations where
    all holomorphic functions are represented by power series of a
    linear polynomial solution)? 
  \item In section \ref{sec_find}, the problem of finding systems
    equipped with $*-$multiplication was treated. But the opposite problem
    is also interesting. Given a linear system of PDE's, determine
    if it has a multiplicative structure on the space of solutions.
  \item Theorem \ref{mult_thm} characterizes all systems that admit
    $*-$multiplication. We have also presented some typical (generic)
    systems with multiplication. However, in order to gain a better
    understanding of the systems that admit multiplication, it would be
    desirable to have some natural principle of classifying these systems. 
  \item In example \ref{jodeitex}, we saw that a particular matrix
    equation of the form $\nabla f = M\nabla g$ admits a
    $*-$multiplication and that the general analytic solution can be
    expressed by power series, with respect to the $*-$multiplication, of
    simple solutions. This is a much more general property, that can be
    generalized to every equation $\nabla f = M\nabla g$ where $M$ is a
    constant matrix with either real or complex entries. The work on this
    problem is already in progress, and the results are being prepared for
    publication.
\end{enumerate}

\section*{Acknowledgments}

I would like to thank Prof. Stefan Rauch-Wojciechowski for comments,
suggestions and helpful discussions.

\raggedright
\bibliography{references}
\end{document}